\documentclass[12pt,letterpaper]{amsart}
\usepackage{palatino, euler, epic, eepic, amssymb, xypic, microtype, enumerate, stmaryrd, epsfig}
\usepackage[outdir=./]{epstopdf}

\usepackage{graphicx,float}
\input xy
\xyoption{all}

 \newlength{\baseunit}               
 \newcount{\numlines}                
\usepackage{amsmath, amssymb, amscd, verbatim, xspace,amsthm}
\usepackage{latexsym, color}

\newcommand{\Z}{\ensuremath{{\mathbb{Z}}}\xspace}
\renewcommand{\P}{\ensuremath{{\mathbb{P}}}}

\newcommand{\R}{\ensuremath{{\mathbb{R}}}}

\newcommand{\lra}{\longrightarrow}

\newcommand\isom{\cong}

\newcommand\Spec{\operatorname{Spec}}

\newcommand\bq{\begin{equation}}
\newcommand\eq{\end{equation}}

\newtheorem{proposition}{Proposition}[section]
\newtheorem{theorem}[proposition]{Theorem}

\newtheorem{example}[proposition]{Example}

\newtheorem{lemma}[proposition]{Lemma}

\theoremstyle{definition}

\theoremstyle{remark}
\newtheorem{remark}[proposition]{Remark}
\usepackage{url}

\numberwithin{equation}{section}

\newcommand{\cut}[1]{}

\newcommand\hidden[1]{}

\newcommand{\FF}{\mathbb{F}}

\newcommand{\PP}{\mathbb{P}}
\newcommand{\QQ}{\mathbb{Q}}
\newcommand{\RR}{\mathbb{R}}
\renewcommand{\R}{\mathbb{R}}
\renewcommand{\Z}{\mathbb{Z}}

\renewcommand{\P}{\mathbb{P}}

\newcommand{\ZZ}{{\mathbb{Z}}}                                        %
\newcommand{\cO}{{\mathcal O}}                                        %
\newcommand{\Cox}{\operatorname{Cox}}                                   %
\newcommand{\Bl}{\operatorname{Bl}}                                   
\newcommand{\Cl}{\operatorname{Cl}}                                   
\usepackage{caption}                                                  %
\usepackage{subcaption}                                               %
\usepackage{tikz-cd}                                                  %
\usetikzlibrary{calc}                                                 %
\usetikzlibrary{positioning}                                          %
\usetikzlibrary{decorations.pathreplacing}                            %
\usetikzlibrary{angles}                                               %
\usetikzlibrary{quotes}                                               %
\usetikzlibrary{shapes.misc}                                          %
\usetikzlibrary{decorations.text}                                     %
\usetikzlibrary{arrows}                                               %
\usetikzlibrary{spy}                                                  %
\title{On a family of negative curves.}

\author{Javier Gonz\'alez Anaya, Jos\'e Luis Gonz\'alez and Kalle Karu}
\address{
J. Gonz\'alez-Anaya, Department of Mathematics, University of British Columbia, 
  Vancouver, BC V6T1Z2, Canada.  \newline \indent
J.L. Gonz\'alez,  Department of Mathematics, University of California, Riverside,
  Riverside, CA 92521, United States.  \newline \indent
K. Karu,
Department of Mathematics, University of British Columbia, 
  Vancouver, BC V6T1Z2, Canada.} 
\email{jga@math.ubc.ca, jose.gonzalez@ucr.edu, karu@math.ubc.ca}
\thanks{J. Gonz\'alez-Anaya was supported by CONACyT scholarship 410172. J.L. Gonz\'alez was supported by an NSERC Discovery grant.}

\begin{document}
\begin{abstract}
Let $X$ be the blowup of a weighted projective plane at a general point. We study the problem of finite generation of the Cox ring of $X$. Generalizing examples of Srinivasan and Kurano-Nishida, we consider examples of $X$ that contain a negative curve of the class $H-mE$, where $H$ is the class of a divisor pulled back from the weighted projective plane and $E$ is the class of the exceptional curve. For any $m>0$ we construct examples where the Cox ring is finitely generated and examples where it is not.

\end{abstract}
\maketitle
\setcounter{tocdepth}{1} 




\section{Introduction}

We work over an algebraically closed field $k$ of characteristic zero.

Recall that the Cox ring of a normal projective variety $X$ is (see \cite{HuKeel})
\[ \Cox(X) = \bigoplus_{[D]\in \Cl(X)} H^0(X, \cO_X(D)).\] 
The variety $X$ is a Mori Dream Space (MDS for short) if its Cox ring is a finitely generated $k$-algebra.

In this article we are mainly interested in the case where $X$ is the blowup of a weighted projective plane at a general point,
\[ X= \Bl_{t_0} \PP(a,b,c).\]
Determining for which triples $(a,b,c)$ is $X$ a MDS is a major open problem. It is known for small values of $a,b,c$, for example when $a\leq 4$, $a=6$, or $abc <(a+b+c)^2$ \cite{Srinivasan, Cutkosky}. For larger values of $a,b,c$ we have only sporadic examples of MDS and non-MDS. For example, $\PP(n,n+1,n+2)$ and $\PP(n,n+1,n+3)$ are MDS \cite{Cutkosky}.

Cutkosky in \cite{Cutkosky} proved that the variety $X$ as above is a MDS if and only if 
\begin{enumerate}
\item $X$ contains a negative curve different from the exceptional curve $E$; that is, an irreducible curve $C\neq E$ such that $C\cdot C\leq 0$.
\item $X$ contains a nonempty curve $D$ disjoint from $C$.
\end{enumerate}
A similar condition for $X$ to be a MDS was also proved by Huneke \cite{Huneke} in terms of the homogeneous coordinate ring of $\PP(a,b,c)$. 
Note that we call an irreducible curve $C$ negative even if its self-intersection number is $0$. Such zero curves are included to cover all cases, for example $\Bl_{t_0} \PP^2$ and blowups of toric varieties below.

We do not know any $X$ that does not contain a negative curve $C\neq E$. When constructing examples, we study varieties $X$ that do contain a negative curve and try to establish the existence or non-existence of the curve $D$ disjoint from $C$. Infinite families of non-MDS were first constructed by Goto, Nishida and Watanabe \cite{GNW}. The construction was later generalized in \cite{GK, He}. In all these examples $X$ contains a negative curve $C$ in the class $\pi^* H - E$, where $\pi: X\to \PP(a,b,c)$ is the projection, $H$ is the class of a curve in $\PP(a,b,c)$, and $E$ is the class of the exceptional curve. 

In the present article we are interested in finding examples of MDS and non-MDS among varieties $X$ that contain a negative curve $C$ having the class of the form $\pi^* H-mE$ for $m\geq 1$. The motivation for this work comes from two sets of examples. Srinivasan in \cite{Srinivasan} considered the example $\Bl_{t_0} \PP(5,77,101)$, which is a MDS. Kurano and Nishida in \cite{KuranoNishida} gave a family of examples of non-MDS, the simplest one of them being $\Bl_{t_0} \PP(16,97,683)$. Both of these examples contain a negative curve of class $\pi^* H-mE$ for $m=2$. Our goal is to generalize these examples to arbitrary $m\geq 1$.

We will take the toric point of view and construct the varieties $X$ as follows. Let $T=\Spec k[x^{\pm 1}, y^{\pm 1}]$ be the torus, and let $C^\circ \subset T$ be an irreducible curve that has multiplicity $m$ at the point $t_0=(1,1)$. We compactify $T$ to a toric variety $X_\Delta$ by choosing a rational triangle $\Delta \subset \RR^2$ that contains the Newton polygon of the curve $C^\circ$. Let $C$ be the strict transform of $C^\circ$ in the blowup $X=\Bl_{t_0} X_\Delta$. Then $C$ has class of the form $\pi^* H-mE$. Moreover, $C$ is a negative curve if the area of $\Delta$ is $\leq \frac{m^2}{2}$.

\begin{figure}
\centering 
  \begin{subfigure}[]{0.4\linewidth}

\begin{tikzpicture}[
    scale=1.5,
    axis/.style={ ->, >=stealth'},
    important line/.style={very thick},
    ]

    \draw[axis] (-0.5,0) -- (2.5,0) node(xline)[right] {};    
    \draw[axis] (0,-0.5) -- (0,3.5) node(yline)[above] {};
    \foreach \Point in {
      (0,0),(0,1),(0,2),(0,3),
      (1,0),(1,1),(1,2),(1,3),
      (2,0),(2,1),(2,2),(2,3)
    }
    \draw[fill=black] \Point circle (0.06);

    \coordinate (A) at (0,0);
    \coordinate (B) at (2,3);
    \coordinate (C) at (1.25,-0.2);

    \draw (A) circle (0) node[below left] {$(0,0)$};
    \draw (B) circle (0) node[above right] {$(2,3)$};
    \draw (C) circle (0) node[below right] {$\frac{1}{101}(125,-5)$};

    \draw (A) -- (B) -- (C) -- cycle;

\end{tikzpicture}
    \caption{$\P(5,77,101)$}  
  \end{subfigure}
\qquad
  \begin{subfigure}[c]{0.4\linewidth}

\begin{tikzpicture}[
    scale=1.5,
    axis/.style={ ->, >=stealth'},
    important line/.style={very thick},
    ]

    \draw[axis] (-0.5,0) -- (2.5,0) node(xline)[right] {};    
    \draw[axis] (0,-0.5) -- (0,3.5) node(yline)[above] {};
    \foreach \Point in {
      (0,0),(0,1),(0,2),(0,3),
      (1,0),(1,1),(1,2),(1,3),
      (2,0),(2,1),(2,2),(2,3)
    }
    \draw[fill=black] \Point circle (0.06);

    \coordinate (A) at (-0.2,0) {};
    \coordinate (B) at (2,3) {};
    \coordinate (C) at (1.25,0) {};

    \draw (A) circle (0) node[below left] {$-\frac{1}{16}$};
    \draw (B) circle (0) node[above right] {$(2,3)$};
    \draw (C) circle (0) node[below right] {$1+\frac{25}{97}$};

    \draw (A) -- (B) -- (C) -- cycle;

\end{tikzpicture}
    \caption{$\P(16,97,683)$}  
  \end{subfigure}
\caption{Triangles corresponding to weighted projective planes.}
\label{fig-oldex}
\end{figure}
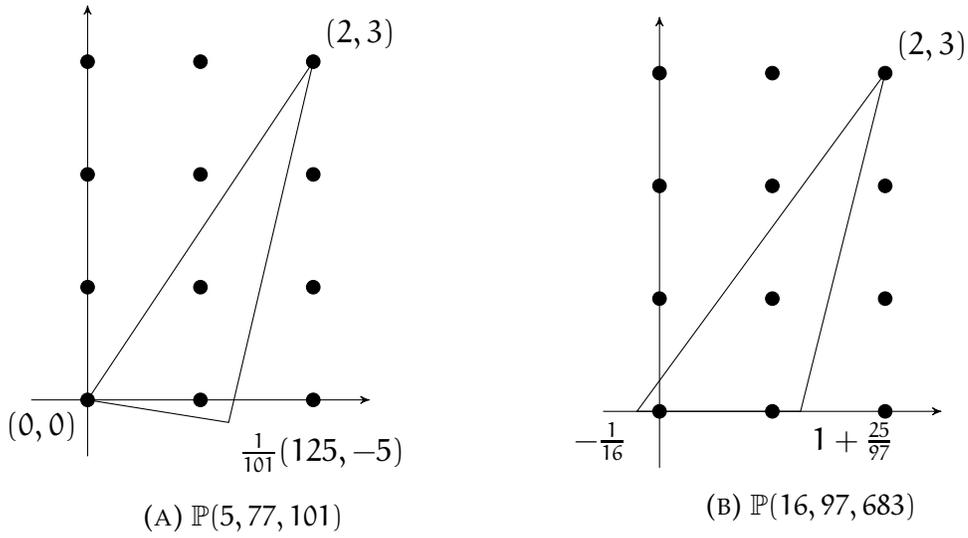

Let us illustrate this construction with the examples of Srinivasan and Kurano-Nishida. In both cases we start with the curve $C^\circ$ defined by the vanishing of the polynomial 
\[ \xi_2(x,y) = 1+x-3xy+x^2y^3.\]
This polynomial is irreducible and has multiplicity $2$ at $t_0$. 
Figure~\ref{fig-oldex} shows the triangles that correspond to $\PP(5,77,101)$ and $\PP(16,97,683)$. Both triangles have area less than $2$, hence $C$ is a negative curve in $X = 
\Bl_{t_0} X_\Delta$. The first example is a MDS, the second one is not.

In order to generalize these examples to arbitrary $m$, we start with the equations defining the curves $C^\circ$.

\begin{theorem} \label{thm-curves}
For any integer $m\geq 1$ there exists an irreducible polynomial, unique up to a constant factor,
\[ \xi_m \in k[x, y] \]
such that
\begin{enumerate}
\item $\xi_m$ vanishes to order $m$ at $t_0=(1,1)$.
\item The Newton polygon of $\xi_m$ is the triangle with vertices $(0,0)$, $(m-1,0)$, $(m,m+1)$.
\end{enumerate}
\end{theorem}

The family of polynomials $\xi_m$ is very regular. Formulas expressing $\xi_{m+1}$ in terms of $\xi_m$ were already given by Kurano and Nishida \cite{KuranoNishida}.

Next we put a triangle $\Delta$ around the Newton polygon of $\xi_m$ and consider $X=\Bl_{t_0} X_\Delta$. Generalizing the example of Srinivasan, we have:
%

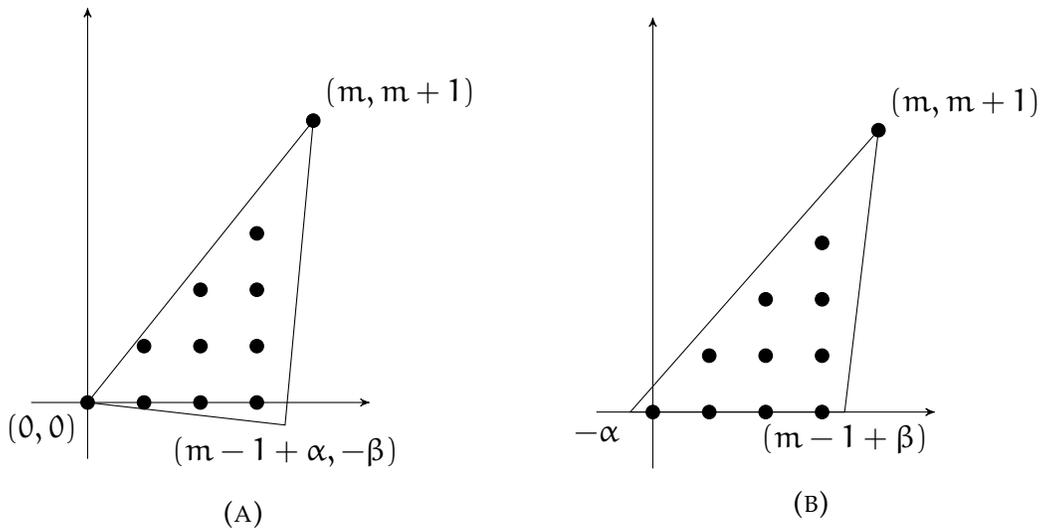
\begin{figure}
\centering 
  \begin{subfigure}[]{0.4\linewidth}

\begin{tikzpicture}[
    scale=1.5,
    axis/.style={ ->, >=stealth'},
    important line/.style={very thick},
    ]

    \draw[axis] (-0.5,0) -- (2.5,0) node(xline)[right] {};    
    \draw[axis] (0,-0.5) -- (0,3.5) node(yline)[above] {};
    \foreach \Point in {
      (0,0),(0.5,0),(1,0),(1.5,0),
      (0.5,0.5),(1,0.5),(1.5,0.5),
      (1,1),(1.5,1),
      (1.5,1.5), (2,2.5)
    }
    \draw[fill=black] \Point circle (0.06);

    \coordinate (A) at (0,0);
    \coordinate (B) at (2,2.5);
    \coordinate (C) at (1.75,-0.2);

    \draw (A) circle (0) node[below left] {$(0,0)$};
    \draw (B) circle (0) node[above right] {$(m,m+1)$};
    \draw (C) circle (0) node[below] {$(m-1+\alpha,-\beta)$};

    \draw (A) -- (B) -- (C) -- cycle;

\end{tikzpicture}
    \caption{}  
  \end{subfigure}
\qquad
  \begin{subfigure}[c]{0.4\linewidth}

\begin{tikzpicture}[
    scale=1.5,
    axis/.style={ ->, >=stealth'},
    important line/.style={very thick},
    ]

    \draw[axis] (-0.5,0) -- (2.5,0) node(xline)[right] {};    
    \draw[axis] (0,-0.5) -- (0,3.5) node(yline)[above] {};
    \foreach \Point in {
      (0,0),(0.5,0),(1,0),(1.5,0),
      (0.5,0.5),(1,0.5),(1.5,0.5),
      (1,1),(1.5,1),
      (1.5,1.5), (2,2.5)
    }
    \draw[fill=black] \Point circle (0.06);

    \coordinate (A) at (-0.2,0);
    \coordinate (B) at (2,2.5);
    \coordinate (C) at (1.7,0);

    \draw (A) circle (0) node[below left] {$-\alpha$};
    \draw (B) circle (0) node[above right] {$(m,m+1)$};
    \draw (C) circle (0) node[below] {$(m-1+\beta)$};

    \draw (A) -- (B) -- (C) -- cycle;

\end{tikzpicture}
    \caption{}  
  \end{subfigure}
\caption{Triangles in Theorems~\ref{thm-main1}-\ref{thm-main2} for $m=4$.}
\label{fig-newex}
\end{figure}  

\begin{theorem} \label{thm-main1}
Let $m\geq 1$ and let $\Delta$ be the triangle (see Figure~\ref{fig-newex}(A)) with vertices 
\[ (0,0), (m-1+\alpha, -\beta), (m,m+1),\]
 where $\alpha,\beta$ are rational numbers such that
\begin{enumerate}
\item The Newton polygon of $\xi_m$ lies in $\Delta$.
\item The area of $\Delta$ is $\leq \frac{m^2}{2}$.
\end{enumerate}
Then $\xi_m$ defines a negative curve $C$ in $X=\Bl_{t_0} X_\Delta$ with class of the form $\pi^* H-mE$. The variety $X$ is a MDS.
\end{theorem}

Generalizing the examples of Kurano and Nishida:

\begin{theorem} \label{thm-main2}
Let $m\geq 1$ and let $\Delta$ be the triangle (see Figure~\ref{fig-newex}(B)) with vertices \[ (-\alpha,0), (m-1+\beta, 0), (m,m+1),\]
 where $\alpha,\beta$ are rational numbers satisfying
\begin{align*} 0 < &\alpha < \frac{1}{1+ (m+1)+(m+1)^2},\\
 \frac{1}{m+2} < &\beta < 1-\frac{1}{1+ \frac{1}{m+1}+\frac{1}{(m+1)^2} - \frac{\alpha}{1-(m+2)\alpha}}.
 \end{align*}
Then $\xi_m$ defines a negative curve $C$ in $X=\Bl_{t_0} X_\Delta$ with class of the form $\pi^* H-mE$. The variety $X$ is not a MDS.
\end{theorem}
A rational triangle $\Delta$ defines a projective toric surface of Picard number $1$. It is isomorphic to some $\PP(a,b,c)$ if the primitive lattice point generators of the rays in the normal fan of $\Delta$ also generate $\ZZ^2$. 

\begin{example}
In Theorem~\ref{thm-main2} the numbers
\[ \alpha=\frac{1}{(m+2)^2}, \quad \beta=\frac{(m+2)^2+1}{(m+2)^3+1}\]
satisfy the inequalities for any $m\geq 1$. The normal fan of the polygon $\Delta$ has rays generated by 
\[ ((m+2)^3+1, -(m+2)^2), \quad (-(m+2)^2, m^2+3m+1), \quad (0,-1).\]
Then 
\[ X_\Delta = \PP((m+2)^2,(m+2)^3+1,(m+2)^3(m^2+2m-1)+m^2+3m+1),\] 
and its blowup is not a MDS. For $m=1,2,3,4,5$ this gives examples
\[ \PP(9,28,59), \quad \PP(16, 65, 459), \quad \PP(25, 126, 1769), \quad \PP( 36, 217,  4997), \quad \PP(49,344, 11703).\]
\end{example}

\begin{remark}
There are toric varieties $X_\Delta$ other than those in Theorems~\ref{thm-main1}-\ref{thm-main2} for which $\xi_m$ defines a negative curve in the blowup. For example, when the triangle $\Delta$ is as in Figure~\ref{fig-newex}(B) with $0\leq \beta\leq \frac{1}{m+2}$ and $\alpha$ small enough so that the triangle has area $\leq \frac{m^2}{2}$, then $X=\Bl_{t_0} X_\Delta$ is a MDS. Indeed, $\xi_{m+1}$ defines a curve $D$ in $X$ that is disjoint from $C$.
\end{remark}

{\bf Acknowledgment.} We thank Ana-Maria Castravet and Zhuang He for telling us about the example by Srinivasan.

\section{Preliminaries about toric varieties and Cox rings.}

Let $X$ be a normal surface. We denote by $\Cl(X)$ the class group of $X$ (Weil divisors modulo linear equivalence) and by $N_1(X)$ the $\RR$-vector space of numerical equivalence classes of curves in $X$. Every curve in $X$ has a class in $\Cl(X)$ and its image in $N_1(X)$. We denote the intersection product of curves by $C\cdot D$. This product is defined between classes in $N_1(X)$.

\subsection{Weil divisors and rational triangles.}

Let $X$ be a normal variety and $D$ an integral Weil divisor in $X$. Recall that the sheaf $\cO_X(D)$ is defined so that its sections $s$ on an open set $U\subseteq X$ are rational functions $f\in K(X)$ such that $div(f)+D$ is effective on $U$. We define the vanishing locus $V(s)$ of the section $s$ as the support of the divisor $div(f)+D$. Notice that the sheaf $\cO_X(D)$ is invertible away from the support of $D$, and more generally, away from $V(s)$ for any global section $s$. If two Weil divisors are linearly equivalent (that means, they have the same class in $\Cl(X)$), then the corresponding sheaves are isomorphic.

Let $X_\Delta$ be a toric surface defined by a rational triangle. Recall that the class group $\Cl(X_\Delta)$ is generated by the torus-invariant divisors. If a toric Weil divisor $D$ is ample, then it corresponds to a rational triangle $\Delta_D$, with sides parallel to the sides of $\Delta$ and each side satisfying an integrality condition. A side is integral if the line that contains the side also contains a lattice point. Two such ample Weil divisors are linearly equivalent if the corresponding triangles differ by an integral translation. The class group may have torsion, as two triangles  may differ by a rational translation.

Let $C^\circ \subset T$ be a curve and $\overline{C}\subset X_\Delta$ its closure. Then $C$ defines a Weil divisor in $X_\Delta$ and its class in $\Cl(X_\Delta)$ is the smallest triangle (with sides parallel to the sides of $\Delta$) that contains the Newton polygon of an equation defining $C^\circ$. Such a triangle has all sides integral, as every side clearly contains a lattice point.

Let $D$ be a divisor with class corresponding to a triangle $\Delta'$. Then the space of global sections of $\cO_{X_\Delta}(D)$ is isomorphic to the space of Laurent polynomials with Newton polygon in $\Delta'$. Let $P\in X_\Delta$ be a $T$-fixed point, corresponding to a vertex of the triangle $\Delta$. The divisor given by a global section of $\cO_{X_\Delta}(D)$ does not pass through the point $P$ if and only if the corresponding vertex of $\Delta'$ is integral and the monomial corresponding to this vertex occurs with nonzero coefficient in the global section.

\subsection{Cones of curves}

Let $X=\Bl_{t_0} X_\Delta$, where $X_\Delta$ is the toric variety defined by a rational triangle. We consider the dual cones $\overline{NE(X)}$ and $Nef(X)$ in $N_1(X) \isom \RR^2$.

Let us generalize the result of Cutkosky mentioned in the introduction from $\PP(a,b,c)$ to $X_\Delta$. The proof is the same as in \cite{Cutkosky}.

\begin{lemma}  \label{lem-Cutkosky}
$X$ is a MDS if and only if
\begin{enumerate}
\item $X$ contains a negative curve $C$ different from the exceptional curve $E$.
\item $X$ contains a nonempty curve $D$ disjoint from $C$.
\end{enumerate}
\end{lemma}

\begin{proof}
The cone $\overline{NE(X)}$ has extremal rays $E \RR_{\geq 0}$ and $\gamma \RR_{\geq 0}$ for some class $\gamma=\pi^* H-cE$, $c>0$. Here $\gamma\cdot\gamma\leq 0$ because otherwise $\gamma$ would satisfy Kleiman's criterion for ampleness.

Assume that $X$ is a MDS. Then $NE(X) = \overline{NE(X)}$, hence the ray  $\gamma \R_{\geq0}$ is generated by the class of a curve $C$, $C\cdot C\leq 0$. We may replace $C$ with a component of $C$ and assume that $C$ is irreducible. This shows that $X$ contains a negative curve $C$ different from $E$.

Consider the extremal ray of $Nef(X)$, $\delta\RR_{\geq 0}$ such that $\delta\cdot\gamma=0$. Since $X$ is a MDS, every nef divisor is semiample. By the lemma below, semiampleness of $\delta$ is equivalent to the existence of a curve $D$, with class on the ray $\delta\RR_{\geq 0}$, such that $D$ is disjoint from $C$.

Conversely, assume that the two conditions are satisfied. Again by the lemma below, $D$ is semiample. Lemma~3 in \cite{Cutkosky} states that if $H$ and $D$ are semiample on a surface $X$, then the ring
\[ \oplus_{m,n\geq 0} \cO_X(mH+nD) \]
is finitely generated. This implies that the Veronese subring of $\Cox(X)$, with degrees in $Nef(X)$, is finitely generated. Any other effective divisor $F$ intersects either $E$ or $C$ negatively, hence it contains $E$ or $C$, and so either $F-E$ or $F-C$ is effective. Thus, adding the defining equations of $E$ and $C$ to the generators of the Veronese subring gives the generators of $\Cox(X)$. 
\end{proof}

\begin{lemma}
Let $C$ be a negative curve in $X$ with class
\[ [C] = \pi^*H - mE \in N_1(X).\]
Consider the class 
\[ \delta = \pi^*H' -mE \in N_1(X)\]
such that $\delta \cdot C = 0$.
Then $\delta$ is semiample if and only if there exists a curve $D$ of class $n\delta$ for some $n>0$, disjoint from $C$. 
\end{lemma}

\begin{proof} 
Assume that $\delta$ is semiample. Then there exists an effective curve $D$ in class $n\delta$ that does not have $C$ as a component. Since $D\cdot C=0$, this $D$ is disjoint from $C$.

Conversely, assume that such a $D$ exists. We need to show that the stable base locus of the divisor $D$ is empty. If $D$ has nonempty stable base locus, then this base locus must be contained in $D$.

The divisor $D-nC= \pi^*(nH'-nH)$ is the pullback of a semiample divisor in $X_\Delta$, and hence is semiample. Since $D-nC$ has empty stable base locus, the stable base locus of $D$ must be contained in $C$. However, $C\cap D$ is empty, hence $D$ is semiample.
\end{proof}

\section{A family of curves.}

We give a proof of Theorem~\ref{thm-curves}, generalizing it slightly. We use notation from \cite{KuranoNishida}.

Consider the following polynomials in $\ZZ[x,y]$:
\[ f=1-xy, \quad g=1-xy^2,\quad h=1-y.\]
Notice that $f,g,h$ vanish at $t_0=(1,1)$.

\begin{proposition}\label{ximRecursion}
There exist polynomials $\xi_m \in\ZZ[x,y]$ for all integers $m\geq 1$ such that
\begin{enumerate}
\item 
\begin{enumerate}
\item $\xi_1 = g$.
\item $\xi_{m+1} = f \xi_m + x^{m} h^{m+1}$.
\item $\xi_{m+1} = x h \xi_m + f^{m+1}$.
\end{enumerate}
\item $\xi_m$ vanishes to order $m$ at $t_0=(1,1)$.
\item The Newton polygon of $\xi_m$ is the triangle with vertices $(0,0)$, $(m-1,0)$, $(m,m+1)$.
\item $\xi_m$ is irreducible in $K[x,y]$, where $K$ is any field.
\end{enumerate}
\end{proposition}

\begin{proof}
Define the polynomials $\xi_m\in\Z[x,y]$ as:
\[
\xi_m=(-1)^mx^my^{m+1} + \sum_{j=0}^{m-1}\sum_{i=j}^{m-1}(-1)^j{m+1\choose j}x^iy^j,
\]
in particular note that $\xi_1=g$ and $\xi_2=1+x-3xy+x^2y^3$. 
The Newton polygon of $\xi_m$ is shown in Figure \ref{xim} below. 

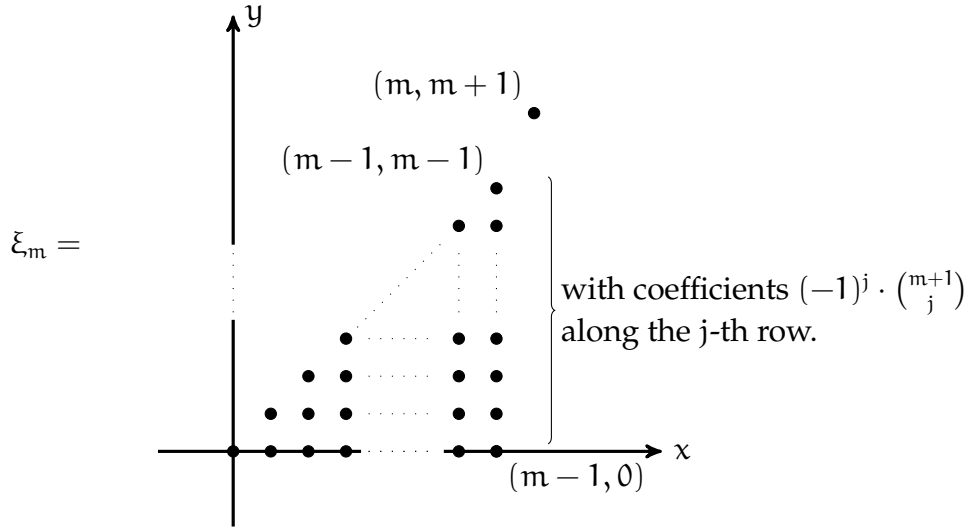
\begin{figure}[h]
\centering

\begin{tikzpicture}[
    scale=1,
    axis/.style={very thick, ->, >=stealth'},
    important line/.style={very thick},
    ]

    \node at (-2,3.25) {$\xi_m=$};
    \draw[important line] (-0.5,0.5) -- (2.2,0.5) node(xline)[right] {};    
    \draw[loosely dotted] (2.3,0.5) -- (3.2,0.5);
    \draw[axis] (3.3,0.5)  -- (6.2,0.5) node(xline)[right] {$x$};

    \draw[important line] (0.5,-0.5) -- (0.5,2.25) node(yline)[above] {};
    \draw[loosely dotted] (0.5,2.35) -- (0.5,3.15);
    \draw[axis] (0.5,3.25)  -- (0.5,6.3) node(xline)[right] {$y$};

    \node at (4.5,4.3) (point) {};
    \draw[decoration={brace,mirror,raise=5.3pt},decorate]
    (4.5,0.6) -- node[right=6pt, align=left] {with coefficients $(-1)^{j}\cdot{m+1\choose j}$\\along the $j$-th row.} (point);

    \foreach \Point/\PointLabel in {
    (0.5,0.5)/{}, (1,0.5)/{},(1.5,0.5)/{}, (2,0.5)/{}, (3.5,0.5)/{},(4,0.5)/{},
    (1,1)/{}, (1.5,1)/{}, (2,1)/{}, (3.5,1)/{},(4,1)/{},
    (1.5,1.5)/{}, (2,1.5)/{}, (3.5,1.5)/{},(4,1.5)/{},
    (2,2)/{}, (3.5,2)/{},(4,2)/{},
    (3.5,3.5)/{},(4,3.5)/{},
    (4,4)/{$(m-1,m-1)$},
    (4.5,5)/{$(m,m+1)$}
    }
    \draw[fill=black] \Point circle (0.075) node[above left] {\PointLabel};

   \foreach \Point/\PointLabel in {(4,0.5)/{$(m-1,0)$} }
    \draw[fill=black] \Point circle (0.075) node[below right] {\PointLabel};
    
    \draw[loosely dotted] (2.3,0.5) -- (3.2,0.5);
    \draw[loosely dotted] (2.3,1) -- (3.2,1);
    \draw[loosely dotted] (2.3,1.5) -- (3.2,1.5);
    \draw[loosely dotted] (2.3,2) -- (3.2,2);
    \draw[loosely dotted] (3.5,2.35) -- (3.5,3.15);
    \draw[loosely dotted] (4,2.35) -- (4,3.15);
    \draw[loosely dotted] (2.15,2.15) -- (3.3,3.3);   
\end{tikzpicture}
\caption{The polynomial $\xi_m$. Monomials in the same row have the same coefficient equal to $(-1)^j\cdot{m+1\choose j}$, where $j$ is the $y$-coordinate of the point.}
\label{xim}
\end{figure}

The proof of (1.b) and (1.c) is then straightforward using the identity ${n\choose k}+{n\choose k-1}={n+1\choose k}$. The algebraic proof of these identities is quite cumbersome and not very helpful, so we present only a diagramatic proof for both in Figure \ref{2identities}.\\

\begin{figure}[H]
\centering 
  \begin{subfigure}[c]{0.4\linewidth}

\begin{tikzpicture}[
    scale=0.4,
    axis/.style={very thick, ->, >=stealth'},
    important line/.style={very thick},
    ]

  \def\s{10}
    \draw[axis] (-1,0)  -- (14,0) node(xline)[right] {$x$};   
    \draw[axis] (0,-1)  -- (0,14) node(xline)[right] {$y$};

    \def\tri#1#2{
      \begin{scope}[shift={#1}]
        \draw[fill= gray, opacity=0.4] (0,0) -- (#2,0) -- (#2,#2) -- cycle;
      \end{scope}
    }
    \def\sq#1{
      \begin{scope}[shift={#1}]
        \def\a{0.3}
        \draw[fill= gray, opacity=0.4] (\a,\a) -- (\a,-\a) -- (-\a,-\a) -- (-\a,\a)-- cycle;
      \end{scope}
    }
    \def\rect#1#2#3{
      \def\a{0.4}
      \begin{scope}[shift={#1}, rotate=#3]
        \draw (-\a,\a) -- (-\a,-\a) -- ($(#2+\a,-\a)$) -- ($(#2+\a,\a)$) -- cycle;
      \end{scope}
      }
    \def\poly#1#2{
      \coordinate (P) at ($#1+(#2+1,#2+2)$);
      \tri{#1}{#2};   
      \sq{(P)};
      }

\poly{(0,0)}{\s}
\poly{(1,1)}{\s}
\def\d{\s+2}
\rect{(\s+1,0)}{\d}{90}
\node at (13.8,7) (point) {$x^mh^{m+1}$};
\node at (7,3) (point) {$f\xi_{m}$};

\end{tikzpicture}
    \caption{$f\xi_m+x^mh^{m+1}$}  
  \end{subfigure}
\qquad
  \begin{subfigure}[c]{0.4\linewidth}

\begin{tikzpicture}[
    scale=0.4,
    axis/.style={very thick, ->, >=stealth'},
    important line/.style={very thick},
    ]
    scale=0.5,
    axis/.style={very thick, ->, >=stealth'},
    important line/.style={very thick},
    ]
  \def\s{10}
    \draw[axis] (-1,0)  -- (14,0) node(xline)[right] {$x$};   
    \draw[axis] (0,-1)  -- (0,14) node(xline)[right] {$y$};

    \def\tri#1#2{
      \begin{scope}[shift={#1}]
        \draw[fill= gray, opacity=0.4] (0,0) -- (#2,0) -- (#2,#2) -- cycle;
      \end{scope}
    }
    \def\sq#1{
      \begin{scope}[shift={#1}]
        \def\a{0.3}
        \draw[fill= gray, opacity=0.4] (\a,\a) -- (\a,-\a) -- (-\a,-\a) -- (-\a,\a)-- cycle;
      \end{scope}
    }
    \def\rect#1#2#3{
      \def\a{0.4}
      \begin{scope}[shift={#1}, rotate=#3]
        \draw (-\a,\a) -- (-\a,-\a) -- ($(#2+\a,-\a)$) -- ($(#2+\a,\a)$) -- cycle;
      \end{scope}
      }
    \def\poly#1#2{
      \coordinate (P) at ($#1+(#2+1,#2+2)$);
      \tri{#1}{#2};   
      \sq{(P)};
      }

\poly{(1,0)}{\s}
\poly{(1,1)}{\s}
\def\d{\s+2}
\rect{(0,0)}{17}{45}

\node at (4.8,7) (point) {$f^{m+1}$};
\node at (7,3) (point) {$xh\xi_{m}$};

\end{tikzpicture}
    \caption{$xh\xi_m+f^{m+1}$}  
  \end{subfigure}
\caption{Visual proofs of (1.b) and (1.c). The highest degree terms of $x^mh^{m+1}$ and $f^{m+1}$ are $(-1)^{m+1}x^{m}y^{m+1}$ and $(-1)^{m+1}x^{m+1}y^{m+1}$, respectively. These terms cancel out the terms $(-1)^mx^my^{m+1}$ of $\xi_m$ and $(-1)^mx^{m+1}y^{m+1}$ of $x\xi_m$, respectively. In both cases the rest of the terms add up precisely to form $\xi_{m+1}$.}
\label{2identities}
\end{figure}
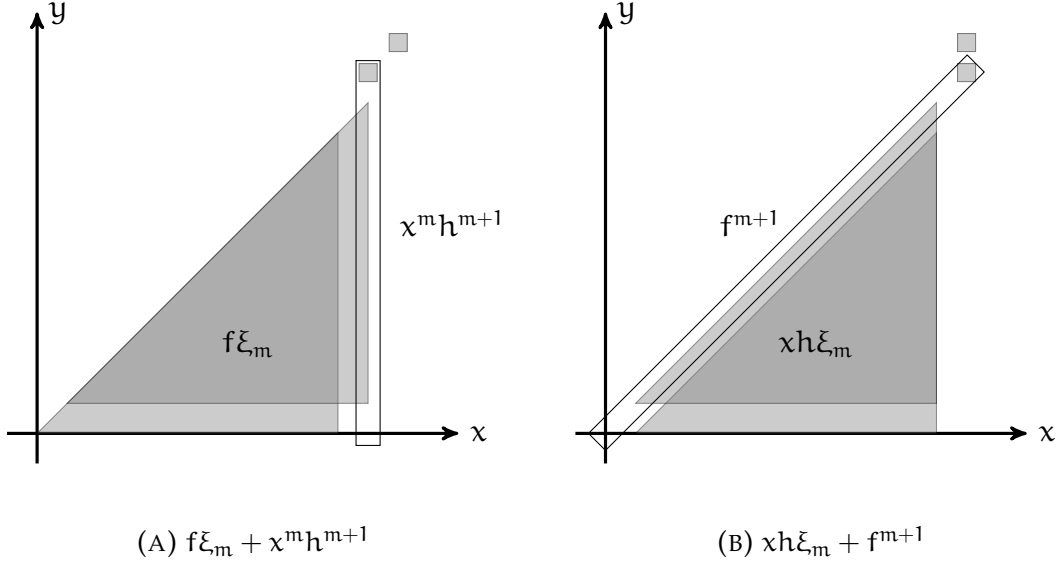  

Proof of (2): To show (2) we argue by induction on $m$. Clearly $f,g$ and $h$ vanish to order $1$ at $t_0$. Indeed, $f=-y(x-1)-(y-1)$ and $g=-y^2(x-1)-(1+y)(1-y)$. Now, assuming $\xi_m$ vanishes to order $m$ at $t_0$, the relation (1.b) (or equivalently (1.c)) implies $\xi_{m+1}$ vanishes to order at least $m+1$ at $t_0$ because both $f\xi_m$ and $x^mh^{m+1}$ vanish to order $m+1$.
The vanishing order of $\xi_{m+1}$ at $t_0$ is exactly $m+1$ since $\partial^{m+1}_x(\xi_{m+1})|_{t_0} =(-1)^{m+1}(m+1)! \neq 0$.

Proof of (3): It follows immediately from the diagrammatic representation of the $\xi_m$.

To prove (4), we first perform a change of coordinates and then apply Eisenstein criterion. Since $\xi_m$ is not divisible by $x$ or $y$, it is irreducible in $K[x,y]$ if and only if it is irreducible in $K[x^{\pm 1}, y^{\pm 1}]$.  
We perform a change of coordinates in $K[x^{\pm 1}, y^{\pm 1}]$ and consider the new polynomial 
\[ \tilde{\xi}_m = x\xi_m(x, y/x).\]
Then the Newton polygon of $\tilde{\xi}_m$ is the triangle with vertices $(1,0)$, $(m,0)$ and $(0,m+1)$. In particular, $\tilde{\xi}_m$ lies in $K[x,y]$ and is irreducible in $K[x,y]$ if and only if it is irreducible in $K[x^{\pm 1}, y^{\pm 1}]$. Now $\tilde{\xi}_m$ can be written as
\[ \tilde{\xi}_m = (-1)^m y^{m+1} + a_m(x) y^m + \cdots + a_0(x),\]
where $x$ divides $a_i(x)$ for $i=1,\ldots,m$ and $x^2$ does not divide $a_0(x)$. Hence $\tilde{\xi}_m$ is irreducible in $K[x,y]$ by Eisenstein criterion.
\end{proof}

As the proof of part (4) of the previous lemma shows, Eisenstein criterion is useful for proving irreducibility of a plane curve if we know the Newton polygon of its equation. 
The following lemma gives a sufficient condition for the irreducibility of a polynomial $f(x,y) \in k[x,y]$ that applies to more general situations.

\begin{lemma} 
Suppose that \[
f(x,y) = \sum_{(i,j) \in P \cap \mathbb{Z}^2}  c_{ij} x^i y^j \in k[x,y],
\] 
where $P$ is a lattice parallelogram $P= \{ w+ \alpha u + \beta v  \ |  \  0 \leq  \alpha \leq m, 0  \leq \beta \leq n   \}$, for some $u,v,w \in \mathbb{Z}^2$ and $m,n \in \mathbb{Z}^{+}$. Then, $f(x,y)$ is irreducible in $k[x,y]$ if the following holds,
\begin{itemize}
\item[(i)] neither $x$ nor $y$ is a factor of f(x,y), 
\item[(ii)] $u,v$ are a basis of $\mathbb{Z}^2$,
\item[(iii)] $c_{w} = c_{w+v} = \cdots = c_{w+(n-1)v} = 0$ and $c_{w+nv+u} = c_{w+nv+2u} = \cdots = c_{w+nv+mu} = 0$. 
\item[(iv)] $c_{w+u} \neq 0$ and $ c_{ w +nv} \neq 0$.
\end{itemize}
\end{lemma}

\section{Proof of Theorem~\ref{thm-main1}.}

Consider the situation in Theorem~\ref{thm-main1}. Let $C$ be the strict transform of the curve $V(\xi_m)$ in $X = \Bl_{t_0} X_\Delta$. Let us first check that $C$ is a negative curve in $X$. The class of $C$ has the form $\pi^*(H)-mE$, where $H$ is the class corresponding to the triangle $\Delta$. Now
\[ C\cdot C = H\cdot H + m^2 E\cdot E  = 2 Area(\Delta) - m^2 \leq 0.\]

By Lemma~\ref{lem-Cutkosky}, it suffices to prove that $X$ contains a curve $D$ disjoint from $C$. We claim that if $D$ is the strict transform in $X$ of the curve $V(1-y)$, then $C\cap D =\emptyset$. 

Let us first check that $C\cdot D=0$.  Let $h$ be the length of the vertical line segment in $\Delta$ from the lower right vertex to the top edge. Since the similar line segment in the smallest triangle containing the Newton polygon of $1-y$ has length $1$, it follows that the class of $D$ in $N_1(X)$ is 
\[ \frac{1}{h}\pi^*(H) - E.\]
Now
\[ C\cdot D = \frac{1}{h} H\cdot H + mE\cdot E = \frac{2}{h} Area(\Delta) - m = 0.\]

Clearly $C$ is not a component of $D$ (as the curves are different when restricted to the torus $T$). This implies that $C$ and $D$ are disjoint.

\section{Characteristic $p$ method of Kurano and Nishida.}

\subsection{Bounding the degree of $D$}

Given a negative curve $C$ in $X$, we are looking for a curve $D$ that is disjoint from $C$. Kurano and Nishida \cite{KuranoNishida} show how to determine a divisor class in which $D$ must lie. We write a weaker version of their result in the toric language.

Consider $X=\Bl_{t_0} X_\Delta$, where $X_\Delta$ is the toric variety given by a rational triangle $\Delta$. Let $C\subset X$ be a negative curve defined by a polynomial $\xi \in\ZZ[x^{\pm 1},y^{\pm 1}]$. Let the class of $C$ be
\[ [C] = \pi^* H -mE \in \Cl(X)\]
for some $H\in\Cl(X_\Delta)$ and integer $m>0$. Assume that $\pi(C)$ passes through a $T$-fixed point $P\in X_\Delta$. (It is not strictly necessary that $P\in \pi(C)$, but it simplifies some proofs below.)

Let $D$ be the class
\[ [D] = \pi^* H' -nE \in \Cl(X)\]
for some $H'\in\Cl(X_\Delta)$ and integer $n>0$ such that $C\cdot D=0$.
 
Let $K$ be a field, either $K=k$ or $K=\overline\FF_p$. We may consider the variety $X$ as well as the polynomial $\xi$ defined over $K$. In particular, $C$ is a negative curve in $X$ that passes through $P$, provided that $char(K)=0$ or $char(K)=p\gg0$.
With notation as above, define
\[ HC_K = \{ l\in \ZZ_{>0} | \cO_X(lD) \text{ has a global section $\zeta$ such that } C \cap V(\zeta) = \emptyset\}.\]
The abbreviation $HC$ stands for Huneke's Condition. The subscript $K$ indicates that we consider $X$ defined over $K$.

We list some properties of $HC_K$ that follow from the definition.
\begin{enumerate}
\item $HC_k$ is nonempty if and only if $X$ is a MDS over $k$.
\item $HC_K$ is closed under addition, hence a sub-semigroup of $(\ZZ_{>0},+)$. This implies that there exists an integer $l_0 \geq 0$ such that $HC_K \subseteq l_0  \ZZ_{>0}$ and $HC_K = l_0  \ZZ_{>0}$ for large numbers.
\item Since $C\cdot D=0$, the condition $C\cap V(\zeta) = \emptyset$ is equivalent to $\zeta$ not vanishing on $C$, which is equivalent to $\zeta$ not vanishing at the $T$-fixed point $P \in C$. Thus, $HC_K$ is defined algebraically as follows. Let $H'$ correspond to a rational triangle $\Delta'$. Then $l\in HC_K$ if and only if there exists a polynomial $\zeta \in K[x^{\pm 1}, y^{\pm 1}]$ with Newton polygon in $l\Delta'$, vanishing to order $ln$ at $t_0$, and with a certain monomial occurring with nonzero coefficient in $\zeta$. The monomial has the exponent given by the vertex of $l\Delta'$ corresponding to the $T$-fixed point $P$.
\item If $l\in HC_K$, then $\cO_X(lD)$ is invertible near $C$ because it has a global section that does not vanish at any point of $C$.
\end{enumerate}

\begin{lemma}
$l\in HC_k$ if and only if $l\in HC_{\overline\FF_p}$ for all $p\gg 0$.
\end{lemma}

\begin{proof}
Let the class $H'$ correspond to a triangle $\Delta'$. Global sections of $\cO_X(lD)$ can then be given as polynomials in $K[x^{\pm 1},y^{\pm 1}]$ that vanish to order $ln$ at $t_0$ and have their Newton polygon contained in $l\Delta'$. The condition of vanishing at $t_0$ is equivalent to the vanishing of partial derivatives up to degree $ln-1$ at $t_0$. This holds when $char(K)=0$ or $char(K)=p \gg 0$ (with $l$ and $n$ fixed). It follows that the space of global sections of $\cO_X(lD)$ is the kernel of an integer matrix $M$, where we view the matrix as a linear map $K^r\to K^s$. The matrix can be put in the Smith normal form (the matrix  becomes diagonal), and then it is clear that its kernel over $K=k$ or $K=\overline\FF_p$ for $p\gg 0$ is equal to the kernel over $\ZZ$ tensored with $K$. In particular, the kernel has the same dimension over any such field.

Similarly, the condition that every global section of $\cO_X(lD)$ vanishes at the $T$-fixed point $P$ is independent of the field $K$. It is equivalent to the condition that all global sections defined over $\ZZ$ vanish at $P$.
\end{proof}

\begin{lemma} 
Let the numbers $m$ and $n$ be as in the definition of $C$ and $D$.
If $l, l+m \in HC_K$ for some $l>0$, then also $m\in HC_K$. 
\end{lemma}

\begin{proof}
Let $\zeta$ be a global section of $\cO_X(lD)$ giving $l\in HC_K$. Consider the exact sequence of sheaves
\[ 0\to \cO_X(mD-nC) \to  \cO_X((l+m)D-nC)\oplus \cO_X(mD) \stackrel{(\xi^n,\zeta)}{\lra} \cO_X((l+m)D) \to 0.\]
Here we fix the sheaf $\cO_X((l+m)D)$ and consider the sheaves $\cO((l+m)D-nC)$ and  $\cO(mD)$ as subsheaves of sections vanishing to order $n$ along $C$ and along the divisor of $\zeta$, respectively. Since the divisor of $\zeta$ does not intersect $C$, the exactness of the sequence is clear.

We claim that the sequence stays exact after applying the global section functor: 
\[ H^0(\cO_X((l+m)D)) = \xi^n H^0(\cO_X((l+m)D-nC))+ \zeta H^0(\cO_X(mD)).\]
Let us prove that this claim implies the statement of the lemma. Indeed, if $\gamma$ is a global section of $\cO_X((l+m)D)$ giving $l+m\in HC_K$, then we can write
\[ \gamma = \xi^n f + \zeta g \]
for some  $f \in H^0(\cO_X((l+m)D-nC))$ and $g\in H^0(\cO_X(mD))$. Let us check that $g$ does not vanish at the point $P\in C$, hence giving $m\in HC_K$. From $\gamma(P)\neq 0$ and $\xi^n(P) = 0$ we get $g(P)\neq 0$. (Note that $\cO_X((l+m)D)$ is invertible near $P$ because it has a global section $\gamma$ that does not vanish at $P$, hence it makes sense to evaluate sections of $\cO_X((l+m)D)$ at $P$.)

To prove the equality of the spaces of global sections above, it suffices to show that 
\[ H^1(\cO_X(mD-nC)) = 0.\]
Here
\[ mD - nC = \pi^*(mH' - nH), \]
where $A=mH' - nH$ is nef on $X_\Delta$ (and even ample if $C\cdot C <0$). This follows from the fact that $D$ lies in the cone of effective curves generated by $[C]$ and $E$. Now
\[ H^1(\cO_X(mD-nC)) =  H^1(\cO_{X_\Delta}(A)) = 0,\]
where the vanishing of $H^1(\cO_{X_\Delta}(A))$ for a nef $\QQ$-Cartier divisor $A$ on the toric variety $X_\Delta$ follows by the same argument as for nef Cartier divisors \cite{Fulton}.
\end{proof}

\begin{proposition} \label{prop-bound}
Assume that $p\in HC_{\overline{\FF}_p}$ for all $p\gg 0$. Then $HC_k$ is nonempty if and only if $m\in HC_k$.
\end{proposition}

\begin{proof}
Assume that $HC_k$ is nonempty, $l_0\in HC_k$. Then also $l_0\in HC_{\overline\FF_p}$ for $p\gg 0$. We know that $HC_{\overline\FF_p} \subseteq l_p \ZZ$ for some $l_p$, and equality holds for large numbers. From $l_0, p\in HC_{\overline\FF_p}$ we get that $l_p| l_0$ and $l_p|p$, hence $l_p=1$ for $p\gg 0$.

Fix a $p\gg 0$. Since $HC_{\overline\FF_p} = \ZZ$ for large numbers, $l, l+m \in HC_{\overline\FF_p}$ for $l$ large. The previous lemma now implies that $m\in HC_{\overline\FF_p}$.

Since $m\in HC_{\overline\FF_p}$ for all $p\gg 0$, it follows that $m\in HC_k$.
\end{proof}

\begin{remark}
Kurano and Nishida prove a stronger statement: with the assumption as in the proposition and $m=1$ or $2$, 
$HC_k$ is nonempty if and only if $1\in HC_k$. In the situation considered below, where $n=m+1$,  this stronger statement for any $m$ would follow from Kawamata-Viehweg vanishing 
\[ H^1(\cO_X(\pi^*(A)) \otimes \omega_X) = 0,\]
where $A$ is an ample divisor in $X_\Delta$. We need this vanishing theorem on the singular variety $X$ in characteristic $p$. The vanishing theorem holds for toric varieties in any characteristic. However, here $X$ is a blowup of a toric variety.
\end{remark}

\subsection{Proof of Theorem~\ref{thm-main2}}

We will apply Proposition~\ref{prop-bound} to prove Theorem~\ref{thm-main2}. The steps in the proof are the same as in \cite{KuranoNishida}.

Let $C$ be the curve in $X$ defined by the polynomial $\xi_m$. Then $C$ has class
\[ \pi^*(H) -mE,\]
where $H$ corresponds to the triangle $\Delta$.

Let $D$ be the class $\pi^*(H') -(m+1)E$, where $H'$ corresponds to the triangle with sides parallel to the sides of $\Delta$ and with base the interval $[0,m]$ on the $x$-axis. 

Let $P$ be the torus fixed point corresponding to the lower left vertex of the triangle $\Delta$. Note that the curve $C$ passes through $P$. A curve of class $lD$ does not pass through $P$ if the constant coefficient of the equation of the curve is nonzero.

\begin{lemma}
$C$ is a negative curve in $X$, $C\cdot C <0$, and $C\cdot D=0$.
\end{lemma}

\begin{proof}
To prove that $C$ is a negative curve, we need to show that the triangle $\Delta$ has area $< \frac{m^2}{2}$, equivalently that $\alpha+\beta < 1/(m+1)$. This can be checked by a direct computation from the inequalities satisfied by $\alpha$ and $\beta$. We set $\beta=\beta(\alpha)$ equal to the upper bound for $\beta$ and find the maximum of the function $\alpha+\beta(\alpha)$. There is a unique critical point on the given interval for $\alpha$, namely at
\[ \alpha = \frac{m+2}{(m+1)^2+(m+2)(1+(m+1)(m+2)^2)}.\]
One now checks that $\alpha+\beta < 1/(m+1)$  at the critical point and at the two endpoints of the interval. 

Comparing the triangles defining $H$ and $H'$, we see that in $N_1(X)$
\[ D = \frac{m}{m-1+\alpha+\beta} H - (m+1)E.\]
Hence
\[ C\cdot D = \frac{m}{m-1+\alpha+\beta} H\cdot H + m(m+1)E\cdot E = \frac{2m}{m-1+\alpha+\beta} Area(\Delta) - m(m+1) = 0.\]
\end{proof}

In the proofs below we will consider triangles with sides parallel to the sides of $\Delta$ and base an interval $[a,b]$ on the $x$-axis. Let us say that a Laurent polynomial $p(x,y)$ lies in degree $[a,b]$ if its Newton polygon lies in such a triangle. As an example, the polynomial $\xi_m$ lies in degree $[-\alpha, m-1+\beta]$. Similarly, $HC_K$ contains an integer $l>0$ if and only if there exists a polynomial $\zeta(x,y)$ that lies in degree $[0,lm]$, vanishes to order $l(m+1)$ at $t_0$, and has nonzero constant term.

To prove Theorem~\ref{thm-main2} using Proposition~\ref{prop-bound}, we show that $p\in HC_{\overline{\FF}_p}$ for all prime numbers $p\gg 0$, but $m\notin HC_k$.

\begin{lemma} 
Let $\Delta$ be as in Theorem~\ref{thm-main2} and $C, D$ as described above. Then $p\in HC_{\overline\FF_p}$ for all $p\gg 0$.
\end{lemma}

\begin{proof}
For any $p\gg 0$ we need to find a polynomial $\zeta=\zeta_p$ that lies in degree $[0,pm]$, vanishes to order $p(m+1)$ at $t_0$, and has nonzero constant term. To construct it, begin by considering the relations (1.b) and (1.c) in Proposition \ref{ximRecursion}. Write $p=(m+1)k+l$, for some $k\in\Z_{\geq 0}$ and $0\leq l \leq m$. Then
\[ (\xi_{m+1}-x^{m}h^{m+1})^p=(f\xi_m)^p,\]
\[
\xi_{m+1}^p+(-x^{m}h^{m+1})^p=(f\xi_m)^{(m+1)k+l}=(f^{m+1})^kf^l\xi_m^p=(\xi_{m+1}-xh\xi_m)^kf^l\xi_m^p,
\]
so that
\begin{align}\label{eqn1}
\xi_{m+1}^p+(-x^{m}h^{m+1})^p=\sum_{i=0}^k{k\choose i}(-1)^i\xi_{m+1}^{k-i}(xh\xi_m)^i\cdot f^l\xi_m^p.
\end{align}
Then $\zeta_p$ is constructed by redistributing the terms in equation \ref{eqn1}:
\begin{align*}
\zeta_p:&=\xi_{m+1}^p-\sum_{i=j+1}^k{k\choose i}(-1)^i\xi_{m+1}^{k-i}(xh\xi_m)^i\cdot f^l\xi_m^p\\
&=(-1)^{p+1}(x^mh^{m+1})^p+\sum_{i=0}^j{k\choose i}(-1)^i\xi_{m+1}^{k-i}(xh\xi_m)^i\cdot f^l\xi_m^p.
\end{align*}
We claim that for $p\gg 0$, there exists $0\leq j \leq k$ such that the left hand side of the equation lies in degree $[0, B]$ and the right hand side lies in degree $[-A,pm]$ for some $A,B>0$. Thus, $\zeta_p$ lies in the intersection of these degrees, $[0,pm]$, as required. 

Every term in the definition of $\zeta_p$ vanishes to order $p(m+1)$ at $t_0$, hence so does $\zeta_p$. Moreover, since $0\leq j$, the constant term of $\zeta_p = \xi_{m+1}^p-\sum(\cdots)$ is equal to the constant term of $\xi_{m+1}^p$, which is nonzero. This proves that $p\in HC_{\overline\FF_p}$.


\begin{figure}
  \centering
  \begin{subfigure}[b]{1\textwidth}
    \centering

\begin{tikzpicture}[
    scale=1,
    axis/.style={ ->, >=stealth'},
    important line/.style={very thick},
    spy using outlines={circle, magnification=3, size=1cm, connect spies}
    ]


    \begin{scope}[shift={(-6,0.5)}]
      \draw[axis] (-0.5,0)  -- (3,0) node(xline)[right] {};   
      \draw[axis] (0,-0.5)  -- (0,3) node(xline)[right] {};
      \foreach \Point in {
        (0,0),(0.5,0),(1,0),(1.48,0),
        (0.5,0.5),(1,0.5),(1.5,0.5),
        (1,1),(1.5,1),
        (1.5,1.5),
        (2,2.5)
      }
      \draw[fill=black] \Point circle (0.05);
      \draw (0,0) -- (1.50125,0) -- (2.04235,2.71235) -- cycle;
      \node at (1,3.75) (point) {$\xi_{m+1}$};
      \node at (1.75,-1) (point) {$\left[0,m+\frac{(m+2)\beta-1}{m+1}\right]$};
    \end{scope}

    \begin{scope}[shift={(-1,0.5)}]
      \draw[axis] (-0.5,0)  -- (3,0) node(xline)[right] {};   
      \draw[axis] (0,-0.5)  -- (0,3) node(xline)[right] {};
      \foreach \Point in {
        (0,0),(0.5,0),(1,0),
        (0.5,0.5),(1,0.5),
        (1,1),
        (1.5,2)
      }
      \draw[fill=black] \Point circle (0.05);
      \begin{scope}[shift={(0.5,0.5)}]]
      \foreach \Point in {
        (0,0),(0.5,0),(1,0),
        (0.5,0.5),(1,0.5),
        (1,1),
        (1.5,2)
      }
      \draw[fill=black] \Point circle (0.05);
      \end{scope}
      \draw (-1/44,0) -- (1.50125,0) -- (2.04065,2.71013) -- cycle;
      \node at (1,3.75) (point) {$f\xi_m$};
      \node at (1.75,-1) (point) {$\left[-\alpha,m+\frac{\beta(m+2)-1}{m+1}\right]$};
    \end{scope}

    \begin{scope}[shift={(3.75,0.5)}]
      \draw[axis] (-0.5,0)  -- (3,0) node(xline)[right] {};   
      \draw[axis] (0,-0.5)  -- (0,3) node(xline)[right] {};
      \begin{scope}[shift={(1.5,0)}]
      \foreach \Point in {
        (0,0),(0,0.5),(0,1),(0,1.5),(0,2),
      }
      \draw[fill=black] \Point circle (0.05);
      \end{scope}
      \draw (-1/44,0) -- (1.5,0) -- (2.04065,2.71013) -- cycle;
      \node at (1,3.75) (point) {$x^mh^{m+1}$};
      \node at (1.75,-1) (point) {$\left[-\alpha,m\right]$};
    \end{scope}

    \begin{scope}[shift = {(-4.3,0.5)},opacity=1]
      \draw[fill = white] (0,0) circle (0.4);

      \filldraw[fill=white] 
      (78.7176:0.4) -- (0,0) -- (0:0.4)
      arc[start angle=0, end angle=78.7176, radius=0.4]
      -- cycle;    
      \draw (-0.4,0) -- (0.4,0);
      \draw[fill=black] (-0.2,0) circle (0.08);
    \end{scope}

    \begin{scope}[shift = {(-1.2,0.5)},opacity=1]
      \draw[fill = white] (0,0) circle (0.4);

      \filldraw[fill=white] 
      (52.7157:0.4) -- (0,0) -- (0:0.4)
      arc[start angle=0, end angle=52.7157, radius=0.4]
      -- cycle;    
      \draw (60:0.4) -- (0.2,0) -- (-60:0.4);
      \draw (-0.4,0) -- (0.4,0);
      \draw[fill=black] (0.2,0) circle (0.08);
    \end{scope}

    \begin{scope}[shift = {(3.55,0.5)},opacity=1]
      \draw[fill = white] (0,0) circle (0.4);

      \filldraw[fill=white] 
      (52.7157:0.4) -- (0,0) -- (0:0.4)
      arc[start angle=0, end angle=52.7157, radius=0.4]
      -- cycle;    
      \draw (60:0.4) -- (0.2,0) -- (-60:0.4);
      \draw (-0.4,0) -- (0.4,0);
    \end{scope}


     \node at (-2.5,2) (point) {$=$};
     \node at (2.5,2) (point) {$+$};
     \node at (-2.5,4.25) (point) {$=$};
     \node at (2.5,4.25) (point) {$+$};

     \begin{scope}[shift={(0,-4.2)}]
       \draw[|-|] (-5,2) -- (-3.5,2);
       \draw (-3.5,2) -- (3.5,2);
       \draw[|-|] (3.5,2) -- (5,2);
       \draw[|-|] (-3.5,1.2) -- (5,1.2);
       \draw[|-|] (-5,0.4) -- (5,0.4);
       \draw[|-|] (-5,-0.4) -- (3.5,-0.4);
       \draw[dotted] (-3.5,2) -- (-3.5,-0.4);
       \draw[dotted] (3.5,2) -- (3.5,-0.4);
     
     \foreach \Point/\PointLabel in { 
       (-5,2)/{$-\alpha$}, (-3.5,2)/{$0$}, (3.5,2)/{$m$}, (5,2)/{$\qquad\qquad m+\frac{(m+2)\beta-1}{m+1}$}, 
       (0,1.2)/{$\xi_{m+1}$}, 
       (0,0.4)/{$f\xi_m$},
       (0,-0.4)/{$x^mh^{m+1}$}, 
     }
     \draw[fill=black] \Point circle (0) node[above] {\PointLabel};
     \end{scope}

\end{tikzpicture}
    \label{}
    \caption{Identity (1.b) in Proposition \ref{ximRecursion} together with the degrees in which all its terms lie.}
  \end{subfigure}
  \begin{subfigure}[b]{1\textwidth}
    \centering

\begin{tikzpicture}[
    scale=1,
    axis/.style={ ->, >=stealth'},
    important line/.style={very thick},
    spy using outlines={circle, magnification=3, size=1cm, connect spies}
    ]


    \begin{scope}[shift={(-6,0.5)}]
      \draw[axis] (-0.5,0)  -- (3,0) node(xline)[right] {};   
      \draw[axis] (0,-0.5)  -- (0,3) node(xline)[right] {};
      \foreach \Point in {
        (0,0),(0.5,0),(1,0),(1.48,0),
        (0.5,0.5),(1,0.5),(1.5,0.5),
        (1,1),(1.5,1),
        (1.5,1.5),
        (2,2.5)
      }
      \draw[fill=black] \Point circle (0.05);
      \draw (0,0) -- (1.50125,0) -- (2.04235,2.71235) -- cycle;
      \node at (1,3.75) (point) {$\xi_{m+1}$};
      \node at (1.75,-1) (point) {$\left[0,m+\frac{(m+2)\beta-1}{m+1}\right]$};
    \end{scope}

    \begin{scope}[shift={(-1,0.5)}]
      \draw[axis] (-0.5,0)  -- (3,0) node(xline)[right] {};   
      \draw[axis] (0,-0.5)  -- (0,3) node(xline)[right] {};
      \begin{scope}[shift={(0.5,0)}]
      \foreach \Point in {
        (0,0),(0.5,0),(1,0),
        (0.5,0.5),(1,0.5),
        (1,1),
        (1.5,2)
      }
      \draw[fill=black] \Point circle (0.05);
      \end{scope}
      \begin{scope}[shift={(0.5,0.5)}]]
      \foreach \Point in {
        (0,0),(0.5,0),(1,0),
        (0.5,0.5),(1,0.5),
        (1,1),
        (1.5,2)
      }
      \draw[fill=black] \Point circle (0.05);
      \end{scope}
      \draw (17/176,0) -- (1601/1000,0) -- (4.2703/2,5.3551/2) -- cycle;
      \node at (1,3.75) (point) {$xh\xi_m$};
      \node at (1.75,-1) (point) {$\left[\frac{1-(m+2)\alpha}{m+1},m+\beta\right]$};
    \end{scope}

    \begin{scope}[shift={(3.75,0.5)}]
      \draw[axis] (-0.5,0)  -- (3,0) node(xline)[right] {};   
      \draw[axis] (0,-0.5)  -- (0,3) node(xline)[right] {};
      \foreach \Point in {
        (0,0),(0.5,0.5),(1,1),(1.5,1.5),(2,2),
      }
      \draw[fill=black] \Point circle (0.05);
      \draw (0,0) -- (1601/1000,0) -- (4.3389/2,5.6989/2) -- cycle;
      \node at (1,3.75) (point) {$f^{m+1}$};
      \node at (1.75,-1) (point) {$\left[0,m+\beta\right]$};
    \end{scope}

    \begin{scope}[shift = {(-4.3,0.5)},opacity=1]
      \draw[fill = white] (0,0) circle (0.4);

      \filldraw[fill=white] 
      (78.7176:0.4) -- (0,0) -- (0:0.4)
      arc[start angle=0, end angle=78.7176, radius=0.4]
      -- cycle;    
      \draw (-0.4,0) -- (0.4,0);
      \draw[fill=black] (-0.2,0) circle (0.08);
    \end{scope}







     \node at (-2.5,2) (point) {$=$};
     \node at (2.5,2) (point) {$+$};
     \node at (-2.5,4.25) (point) {$=$};
     \node at (2.5,4.25) (point) {$+$};

     \begin{scope}[shift={(0,-4.2)}]
       \draw[|-|] (-5,2) -- (-3.5,2);
       \draw (-3.5,2) -- (3.5,2);
       \draw[|-|] (3.5,2) -- (5,2);
       \draw[|-|] (-5,1.2) -- (3.5,1.2);
       \draw[|-|] (-3.5,0.4) -- (5,0.4);
       \draw[|-|] (-5,-0.4) -- (5,-0.4);
       \draw[dotted] (-3.5,2) -- (-3.5,-0.4);
       \draw[dotted] (3.5,2) -- (3.5,-0.4);
     
     \foreach \Point/\PointLabel in { 
       (-5,2.1)/{$0$}, (-3.5,2)/{$\frac{1-(m+2)\alpha}{m+1}$}, (3.5,2)/{$m+\frac{(m+2)\beta-1}{m+1}$}, (5,2)/{$\qquad m+\beta$}, 
       (0,1.2)/{$\xi_{m+1}$}, 
       (0,0.4)/{$xh\xi_m$},
       (0,-0.4)/{$f^{m+1}$}, 
     }
     \draw[fill=black] \Point circle (0) node[above] {\PointLabel};
     \end{scope}

\end{tikzpicture}
    \label{}
    \caption{Identity (1.c) in Proposition \ref{ximRecursion} together with the degrees in which all its terms lie.}
  \end{subfigure}
\caption{Visual representation of the identities in Proposition \ref{ximRecursion} for $m=3$, $\alpha=\frac{1}{22}$ and $\beta=\frac{101}{500}$. The circles represent a zoom in into the area of the figures they cover.}
\label{tri-intervals}
\end{figure}
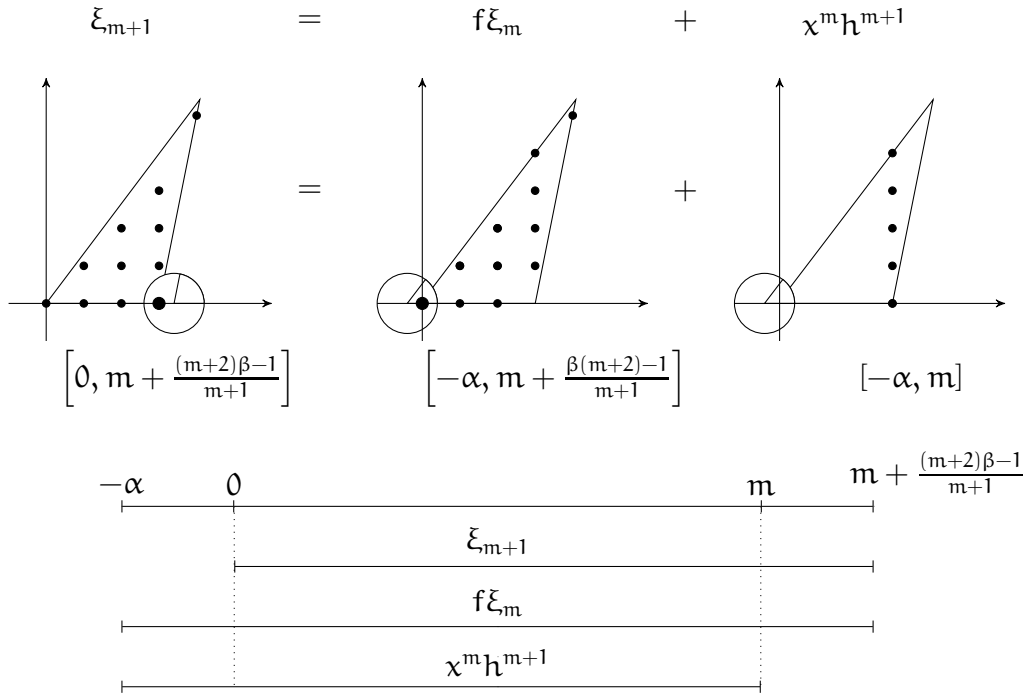
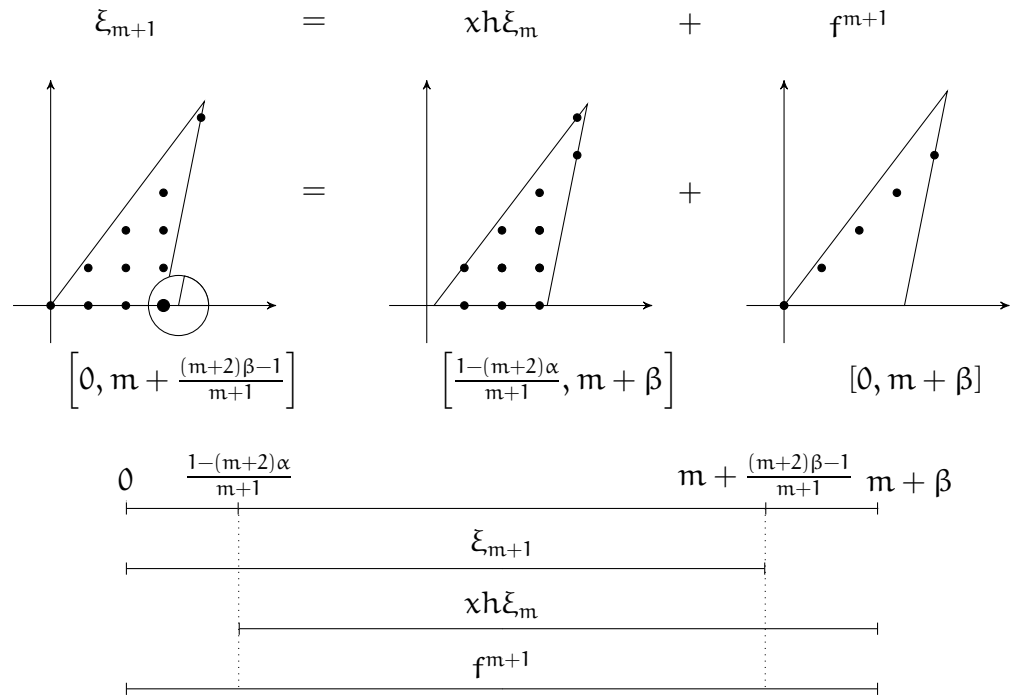

To show the existence of $j$ we begin by studying the degrees in which all the elements in equation (\ref{eqn1}) lie.
Figure \ref{tri-intervals} shows the degree in which each term in the relations (1.b) and (1.c) of Proposition \ref{ximRecursion} lies. We have assumed here that $0<\alpha<\frac{1}{m+2}$ and $\frac{1}{m+2}<\beta<1$, otherwise the degrees would not be as pictured. We will sharpen these bounds as we move forward. 

The terms in the left hand side of equation (\ref{eqn1}) lie in the following degrees:
\begin{align*}
  \xi_{m+1}^p&: p\cdot\left[0,m+\frac{(m+2)\beta -1}{m+1}\right]=\left[0,pm+p\cdot\frac{(m+2)\beta -1}{m+1}\right],\\
(x^{m}h^{m+1})^p&: p\cdot\left[-\alpha,m\right]=\left[-p\alpha,pm\right].
\end{align*}
On the other hand, the degree in which the $i$-th summand of the right hand side of (\ref{eqn1}) lies is determined by the sum of the following degrees:

\begin{align*}
\xi_{m+1}^{k-i}&: (k-i)\cdot\left[0,m+\frac{(m+2)\beta -1}{m+1}\right]=\left[0,(k-i)m+(k-i)\cdot\frac{(m+2)\beta -1}{m+1}\right],\\
(xh\xi_m)^i&: i\cdot\left[\frac{1-(m+2)\alpha}{m+1},m+\beta\right]=\left[i\cdot\frac{1-(m+2)\alpha}{m+1},im+i\beta\right],\\
f^l&:\left[0,l\cdot\frac{m+\beta}{m+1}\right],\\
\xi_m^p&: p\cdot\left[-\alpha, m-1+\beta\right]=\left[-p\alpha,p(m-1)+p\beta\right].
\end{align*}

After some algebra, the $i$-th term in the right hand side of (\ref{eqn1}) lies in degree:
\begin{align*}
\left[-p\alpha+i\cdot\frac{1-(m+2)\alpha}{m+1},pm+p\cdot\frac{(m+2)\beta-1}{m+1}-(k-i)\cdot\frac{1-\beta}{m+1}\right].
\end{align*}

Figure \ref{LHSintervals} summarizes the previous calculations.
\begin{figure}
\centering
\begin{tikzpicture}[
    scale=1,
    ]

\draw[|-|] (-5,2) -- (-3.5,2);
\draw (-3.5,2) -- (3.5,2);
\draw[|-|] (3.5,2) -- (5,2);
\draw[|-|] (-3.5,1) -- (5,1);
\draw[|-|] (-5,0) -- (3.5,0);
\draw[|-|,dashed] (-5,-1.1) -- (-3,-1.1);
\draw (-3,-1.1) -- (3,-1.1);
\draw[|-|,dashed] (3,-1.1) -- (5,-1.1);

\draw[dotted] (-3.5,2) -- (-3.5,0);
\draw[dotted] (3.5,2) -- (3.5,0);

\foreach \Point/\PointLabel in { (0,2)/{},
(-5,2)/{$-p\alpha$}, (-3.5,2)/{$0$}, (3.5,2)/{$pm$}, (5,2)/{$\qquad\qquad pm+p\cdot\frac{\beta(m+2)-1}{m+1}$}, 
(0,1)/{$\xi_{m+1}^p$}, (0,0)/{$(x^mh^{m+1})^p$}, 
(0,-1.1)/{$i$-th term}, (-4,-1.1)/{$i\cdot\frac{1-(m+2)\alpha}{m+1}$}, (4,-1.1)/{$\quad(k-i)\cdot\frac{1-\beta}{m+1}$}
}
\draw[fill=black] \Point circle (0) node[above] {\PointLabel};
\end{tikzpicture}
\caption{All the relevant degrees in Equation \ref{eqn1}.}
\label{LHSintervals}
\end{figure}
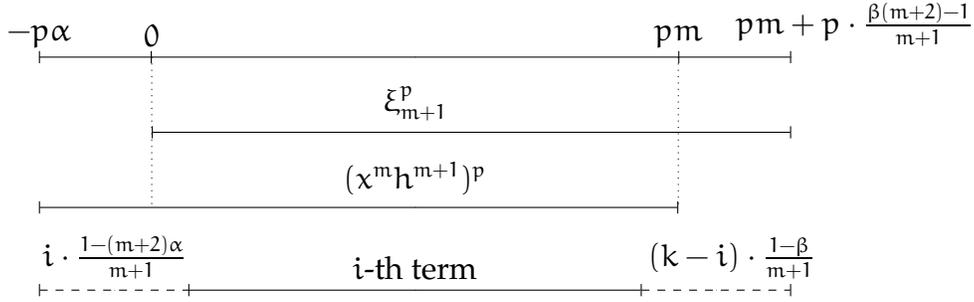

We need to show that for $p\gg 0$ there exists $j\geq 0$ such that the $i$-th term of the sum lies in the same degree as $(x^mh^{m+1})^p$ for $i=0,\dots,j$ and in the same degree as $\xi_{m+1}^p$ for $i=j+1,\dots,k$. This is equivalent to:
\begin{align*}
 p\alpha\leq (j+1)\cdot\frac{1-(m+2)\alpha}{m+1}\qquad\text{and}\qquad p\cdot\frac{(m+2)\beta-1}{m+1}\leq (k-j)\cdot\frac{1-\beta}{m+1}.
\end{align*}
After solving for $j$ and substituting $p=(m+1)k+l$, these inequalities become
\begin{align*}\label{ineq}
\left(\frac{(m+1)^2\alpha}{1-(m+2)\alpha}\right)k+(\cdots)\leq j\leq\left(\frac{(m+2)-((m+1)(m+2)+1)\beta}{1-\beta}\right)k+(\cdots),
\end{align*}
where the terms $(\cdots)$ depend on $l,m,\alpha$ and $\beta$ but do not grow with $k$. The existence of $j$ for large $p$ (that means, large $k$) now follows from the strict inequality  
\[
\frac{(m+1)^2\alpha}{1-(m+2)\alpha}<\frac{(m+2)-((m+1)(m+2)+1)\beta}{1-\beta},
\]
which simplifies to the upper bound for $\beta$ in the statement of Theorem \ref{thm-main2}. Note that the left hand side of the inequality is non-negative, hence we can choose $j\geq 0$. 

The upper bound for $\alpha$,
\[
\alpha=\frac{1}{(m+1)^2+(m+1)+1}<\frac{1}{m+2}.
\]
is the point where the the lower and upper bounds of $\beta$ meet. (For $\alpha$ greater than or equal to this upper bound, there are no values of $\beta$ between the required bounds.)
\end{proof}

\begin{lemma}
Let $\Delta$ be as in Theorem~\ref{thm-main2} and $C, D$ as described above. Then $m\notin HC_k$. 
\end{lemma}
\begin{proof}
Assume $m\in HC_k$, so that there exists a $\zeta\in H^0(\cO(mD))$ that lies in degree $[0,m^2]$, vanishes to order $m(m+1)$ at $t_0$ and has nonzero constant term. Note that $\xi_{m+1}^m$ also vanishes to order $m(m+1)$ at $t_0$ and has nonzero constant term. We show the existence of a triangle $\overline{\Delta}$ such that it contains the Newton polygons of both $\zeta$ and $\xi_{m+1}^m$, and has area less than $\frac{m^2(m+1)^2}{2}$. This implies that  both $\zeta$ and $\xi_{m+1}^m$ define negative curves in $\Bl_{t_0}X_{\overline{\Delta}}$ of the same divisor class. Since $\xi_{m+1}$ is irreducible, $\zeta=c\cdot\xi^{m}_{m+1}$ for some constant $c$. This gives a contradiction as follows.

Let $\Delta_1$ be the triangle of $mH'$ containing the Newton polygon of $\zeta$, and let $\Delta_2$ be the Newton polygon of $\xi_{m+1}^m$. Observe that $\Delta_2$ contains a lattice point at the top vertex that does not lie in $\Delta_1$, see Figure \ref{newtonP}. This contradicts the equality $\zeta = c\cdot\xi^{m}_{m+1}$.

\begin{figure}
\centering

\begin{tikzpicture}[
    scale=2,
    axis/.style={ ->, >=stealth'},
    important line/.style={very thick},
    spy using outlines={circle, magnification=3, size=1cm, connect spies}
    ]



      \draw[thick] (0,0) -- (1.5,0) -- (1.7,2.4) -- cycle;
      \draw[thick] (0,0) -- (1.5,0) -- (2,2.5) -- cycle;
      \node at (0,-0.2) (point) {$0$};
      \node at (1.5,-0.2) (point) {$m^2$};
      \draw (1.7,2.4) circle (0.0) node[above] {$\Delta_1$};
      \draw (2,2.5) circle (0.0) node[above right] {$\Delta_2$};
      \node at (3,1.5) [rectangle] (point1) {$m+2<$};
      \node at (3.6,1.5) [circle] (point2) {$\frac{m+1}{1-\beta}$};
      \def\myshift#1{\raisebox{-2ex}}
      \draw [<-,thick,postaction={decorate,decoration={text along path,text align=center,text={|\myshift|slope}}}] (1.77,1.25) to [bend right =20] (point1);
      \def\myshift#1{\raisebox{1ex}}
      \draw [<-,thick,postaction={decorate,decoration={text along path,text align=center,text={|\myshift|slope}}}] (1.67,12*0.15) to [bend left =30] (point2);
      \node at (0,2) [rectangle] (point3) {$<\frac{m+1}{m+\alpha}$};
      \node at (-0.55,2) [circle] (point4) {$\frac{m+2}{m+1}$};
      \def\myshift#1{\raisebox{1ex}}
      \draw [->,thick,postaction={decorate,decoration={text along path,text align=center,text={|\myshift|slope}}}] (point3) to [bend left =20] (1.25,156/85);
      \def\myshift#1{\raisebox{-2.1ex}}
      \draw [->,thick,postaction={decorate,decoration={text along path,text align=center,text={|\myshift|slope}}}] (point4) to [bend right =10] (1.17,1.5);

\end{tikzpicture}
\caption{The triangle $\Delta_1$ is associated to $\zeta$ and is similar to the triangle defining the original weighted projective space; the triangle $\Delta_2$ is the Newton polygon of $\xi_{m+1}^m$ and its upper vertex is at $m(m+1,m+2)$. Note that given the slopes of the triangles their relative position is as shown in the picture for $\alpha<\frac{1}{m+2}<\beta$.}
\label{newtonP}
\end{figure}
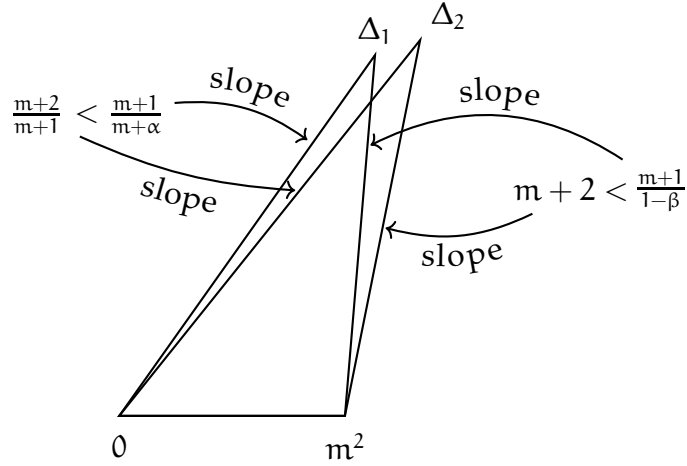
We proceed to show the existence of $\overline{\Delta}$.

Let $\overline{\Delta}$ have base the interval $[0,m^2]$ on the $x$-axis. Take the slope of the right edge of $\overline{\Delta}$ to be the same as that of $\Delta_2$, which is $m+2$. For the slope of the left edge take the line passing through $(0,0)$ and the point $Q=(m^2+1,m(m+1)+1)$. To understand this choice notice that the slope of the left edge of $\Delta_1$ is $\frac{m+1}{m+\alpha}<\frac{m+1}{m}$. Then, the point $Q$ is the first lattice point in $\Delta_1$ the line comes across as we lower the slope starting from $\frac{m+1}{m}$, see Figure \ref{Qconstruction}. It's clear then that the triangle thus obtained contains all lattice points from both $\Delta_1$ and $\Delta_2$.

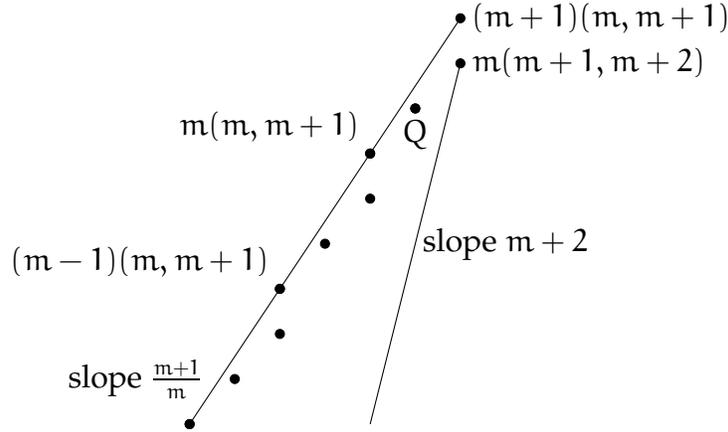
\begin{figure}
\centering

\begin{tikzpicture}[
    scale=0.6,
    axis/.style={ ->, >=stealth'},
    important line/.style={very thick},
    spy using outlines={circle, magnification=3, size=1cm, connect spies}
    ]


      \foreach \Point\PointLabel in {
        (0,0)/{},(1,1)/{},(2,2)/{},
        (2,3)/{},(3,4)/{},(4,5)/{},
        (4,6)/{},(5,7)/{},(6,8)/{m(m+1,m+2)},
        (6,9)/{(m+1)(m,m+1)}
      }
      \draw[fill=black] \Point circle (0.1) node[right] {$\PointLabel$};
      \foreach \Point\PointLabel in {
        (5,7)/{Q},
      }
      \draw[fill=black] \Point circle (0.1) node[below] {$\PointLabel$};
      \foreach \Point\PointLabel in {
        (4,6)/{m(m,m+1)},
        (2,3)/{(m-1)(m,m+1)}
      }
      \draw[fill=black] \Point circle (0.1) node[above left] {$\PointLabel$};
      \draw (0,0) -- (6,9);
      \draw (4,0) -- (6,8);
      \node at (-1.2,1) [rectangle] (point) {slope $\frac{m+1}{m}$};
      \node at (7,4) [rectangle] (point) {slope $m+2$};

\end{tikzpicture}
\caption{The top part of $\overline{\Delta}$ for $m=2$.}
\label{Qconstruction}
\end{figure}

To show that the curves defined by $\xi_{m+1}^m$ and $\zeta$ in $\operatorname{Bl}_{t_0}(X_{\overline{\Delta}})$ are negative we need the area of the triangle $\overline{\Delta}$ to be smaller than $\frac{m^2(m+1)^2}{2}$, i.e., its height must be less than $(m+1)^2$.

The line spanned by the left edge of $\overline{\Delta}$ reaches the height $(m+1)^2$ at the $x$-coordinate $x_L=m(m+1)+\frac{m+1}{m(m+1)+1}$, while the right edge does it at $x_R=m(m+1)+\frac{1}{m+2}$. The height of the triangle is then smaller than $(m+1)^2$ whenever $x_R<x_L$ (see Figure \ref{topDeltaBar}), that is, when
\[
m(m+1)+\frac{1}{m+2}<m(m+1)+\frac{m+1}{m(m+1)+1}\quad\iff\quad \frac{1}{m+2}<\frac{m+1}{m(m+1)+1},
\]
which always holds. 
\begin{figure}
\centering

\begin{tikzpicture}[
    scale=2,
    axis/.style={ ->, >=stealth'},
    important line/.style={very thick},
    spy using outlines={circle, magnification=3, size=1cm, connect spies}
    ]



  \draw (0,0) -- (1.8,2);
  \draw (1,0) -- (1.5,2);
  \draw[dashed] (-0.2,1.8) -- (2,1.8);
  \draw[fill=black] (0.9*0.7,0.7) circle (0.025) node[above left] {$Q$};
  \draw[fill=black] (1.25,1) circle (0.025) node[right] {$m(m+1,m+2)$};
  \draw[fill=black] (1.8*9/10,1.8) circle (0.025) node[below right] {$x_L$};
  \draw[fill=black] (1.5*4/3-0.55,1.8) circle (0.025) node[below left] {$x_R\,$};
  \node at (0.8,0.5) (point) {$\overline{\Delta}$};
  \node at (2.7,1.8) (point) {$y=(m+1)^2$};

\end{tikzpicture}
\caption{The top vertex of $\overline{\Delta}$ always lies below the line $y=(m+1)^2$.}
\label{topDeltaBar}
\end{figure}
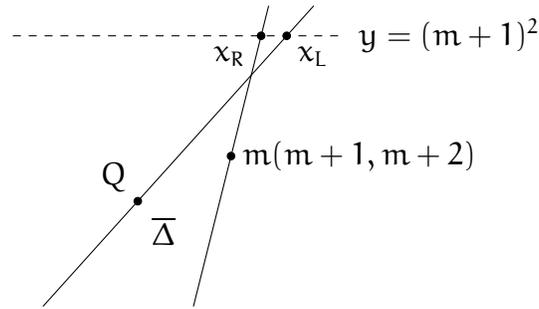
\end{proof}

\begin{remark}
Note from the proofs above that if $\beta\leq \frac{1}{m+2}$ then $\xi_{m+1}$ lies in degree $[0,m]$, and hence defines a curve disjoint from $C$.
\end{remark}


\begin{thebibliography}{1}

\bibitem{Cutkosky}
Steven~Dale Cutkosky.
\newblock Symbolic algebras of monomial primes.
\newblock {\em J. Reine Angew. Math.}, 416:71--89, 1991.

\bibitem{Fulton}
William Fulton.
\newblock {\em Introduction to toric varieties}, volume 131 of {\em Annals of
  Mathematics Studies}.
\newblock Princeton University Press, Princeton, NJ, 1993.
\newblock The William H. Roever Lectures in Geometry.

\bibitem{GK}
Jos\'e~Luis Gonz\'alez and Kalle Karu.
\newblock Some non-finitely generated {C}ox rings.
\newblock {\em Compos. Math.}, 152(5):984--996, 2016.

\bibitem{GNW}
Shiro Goto, Koji Nishida, and Keiichi Watanabe.
\newblock Non-{C}ohen-{M}acaulay symbolic blow-ups for space monomial curves
  and counterexamples to {C}owsik's question.
\newblock {\em Proc. Amer. Math. Soc.}, 120(2):383--392, 1994.

\bibitem{He}
Zhuang He.
\newblock New examples and non-examples of {M}ori {D}ream {S}paces when blowing
  up toric surfaces.
\newblock {\em arXiv: 1703.00819}.

\bibitem{HuKeel}
Yi~Hu and Sean Keel.
\newblock Mori dream spaces and {GIT}.
\newblock {\em Michigan Math. J.}, 48:331--348, 2000.
\newblock Dedicated to William Fulton on the occasion of his 60th birthday.

\bibitem{Huneke}
Craig Huneke.
\newblock Hilbert functions and symbolic powers.
\newblock {\em Michigan Math. J.}, 34(2):293--318, 1987.

\bibitem{KuranoNishida}
Kazuhiko Kurano and Koji Nishida.
\newblock Infinitely generated symbolic {R}ees rings of space monomial curves
  having negative curves.
\newblock {\em arXiv: 1705.09865}.

\bibitem{Srinivasan}
Hema Srinivasan.
\newblock On finite generation of symbolic algebras of monomial primes.
\newblock {\em Comm. Algebra}, 19(9):2557--2564, 1991.

\end{thebibliography}

\end{document}